\documentclass[12pt]{article}

\usepackage[T2A]{fontenc}
\usepackage[cp866]{inputenc}

\usepackage{amsfonts}
\usepackage{amssymb}
\usepackage{amsmath}
\usepackage{amsthm}
\begin{document}

\centerline{\bf A version of the proof for } \centerline{\bf Peres-Schlag's theorem on lacunary sequences.} \vskip+1.5cm \centerline{\bf
Moshchevitin N.G. \footnote{ Research is supported by grants RFFI 06-01-00518, MD-3003.2006.1, NSh-1312.2006.1 and INTAS 03-51-5070
 }
} \vskip+1.5cm

We present a  proof of a multidimensional version of Peres-Schlag's theorem on Diophantine approximations with lacunary sequences. \vskip+0.5cm

 {\bf 1. Introduction.}

A sequence $\{t_j\},\,\,\ j = 1,2,3,..$ of positive real numbers is defined to be lacunary if for
some $ M>0$ one has
\begin{equation}
\frac{t_{j+1}}{t_j} \ge 1+\frac{1}{M},\,\,\,
\forall j \in \mathbb{N}.
\label{lac}
\end{equation}
Erd\"os \cite{E} conjectured that for any lacunary sequence there exists real $\alpha$ such that the set of fractional parts $ \{ \alpha t_j\} ,
\,\, j \in\mathbb{N}$ is not dense in $[0,1]$. This conjecture was proved by Pollington \cite{Po} and de Mathan \cite{M}. Some quantitative
improvements were due to Katznelson \cite{K}, Akhunzhanov and Moshchevitin \cite{AM} and Dubickas \cite{D}. The best known quantitative estimate
is due to Peres and Schlag \cite{P}. The last authors proved that with some positive constant $\gamma > 0$ for any sequence $\{t_j\}$ under
condition (\ref{lac}) there exists a real number $\alpha $ such that
$$
||\alpha t_j  || \ge \frac{\gamma}{M\log M},\,\,\,
\forall j \in \mathbb{N}.
$$
In their proof Peres and Schlag use a special variant of the Lovasz local lemma.

In the present paper we (following the arguments from \cite{P}) prove  a multidi- mensional version of the above result by Peres and Schlag. Our
proofs avoid the terminology from probability theory.

 {\bf 2. Notation and parameters.} $\mu (\cdot )$ denotes the Lebesque measure.
For a set $ A \subset \mathbb{R}^d$ we denote $A^c =[0,1]^d\setminus A.$ Let $M \ge 8$ and $t_1\ge
2$. We shall use the parameters
$$
\delta = \frac{1}{2^{11}(M\log M)^{1/d}} ,\,\,\, \delta_1 = 2^{4d+1}\delta^d, \,\,\,
 h =\lceil2^{3d} M\log M\rceil
$$
and
\begin{equation}
l_j = \left\lfloor \log_2\frac{t_j}{2\delta}\right\rfloor.
\label{DD}
\end{equation}
Note that
\begin{equation}
1-(16\delta )^dh \ge 1/2, \label{de}
\end{equation}
\begin{equation}
 1-{2\delta_1 h} \ge 1/2,
\label{HD}
\end{equation}
and from (\ref{lac}) one has
\begin{equation}
\frac{t_{i+h}}{t_i}\ge \left( 1+\frac{1}{M} \right)^h\ge M^{2^{3d}\log 2} \ge \frac{1}{\delta}. \label{HH}
\end{equation}
 For the proof of our result we need the sets
$$E(j,a) = \left\{ x \in [0,1]:\,\, \left| x - \frac{a}{t_j} \right| \le \frac{\delta}{t_j}\right\}
.$$ For integer point $ {\bf a} = (a_1, . . . , a_d )\in\mathbb{R}^d$ denote
$$E(j,{\bf a}) = E(j,a_1) \times . . . \times  E(j,a_d)$$
 Each $E_{j,a}$ is covered by a  opened dyadic segments of the form
$\left(\frac{b}{2^{l_j}}, \frac{b+\varepsilon}{2^{l_j}}\right),\,\,\,\varepsilon \in \{1,2\}.$ Then
the set
$$
\bigcup_{0\le a_1, . . . , a_d \le \lceil t_j \rceil} E(j,{\bf a}) \bigcap [0,1]^d
$$
is covered by some union of dyadic boxes of the form
$$\left(\frac{b_1}{2^{l_j}},
\frac{b_1+\varepsilon_1}{2^{l_j}}\right)\times . . . \times \left(\frac{b_d}{2^{l_j}},
\frac{b_d+\varepsilon_d}{2^{l_j}}\right) \,\,\,\varepsilon_k \in \{1,2\}.$$
 We denote this union as $A_j$. So
$$
\bigcup_{0\le a_1,...,a_d\le \lceil t_j \rceil} E(j,{\bf a}) \cap [0,1]^d \subseteq A_j,
$$
and  the complement $A_j^c =[0,1]^d\setminus A_j$ can be represented as a union $A_j^c = \cup_{1\le
\nu \le T_j} I_\nu$ of closed dyadic boxes  of the form
\begin{equation}
\left[\frac{b_1}{2^{l_j}}, \frac{b_1+1}{2^{l_j}}\right] \times . . . \times
\left[\frac{b_d}{2^{l_j}}, \frac{b_d+1}{2^{l_j}}\right] .
\label{III}
\end{equation}
 Moreover the set $\cap_{j\le i} A_j^c$
also can be represented as
$$
\bigcap_{j\le i}
A_j^c =\bigcup_{1\le \nu \le T_i} I_\nu
$$
with dyadic intervals $I_\nu$ of the form (\ref{III}).
 Note that
\begin{equation}
\mu (A_j) \le \left( \frac{4\delta}{t_j}\right)^d (\lceil t_j\rceil +1)^d\le (16\delta )^d . \label{mu}
\end{equation}

{\bf 3. The results.}

{\bf Theorem 1.}\,\,\,{\it Let $d \in \mathbb{N}$ and $ t_1\ge 2,\,\, M \ge 8$. Then for any
sequence $\{ t_j\}$ under condition (\ref{lac}) there exists a set of real numbers $
(\alpha_1,...,\alpha_d)$ such that
$$
\max_{1\le  k \le d} ||\alpha_k t_j|| \ge \frac{1}{2^{11}(M\log M)^{1/d}},\,\,\, \forall j \in \mathbb{N}.
$$}
We give the proof of the Theorem 1 in sections 4,5.

For a vector $\xi = (\xi_1,...,\xi_d) \in \mathbb{R}^d$ we define $$ ||\xi || = \max_{1\le j\le d} \min_{a_j \in \mathbb{Z}} |\xi_j - a_j|.$$

{\bf Theorem 2.}\,\,\,{\it Let $d \in \mathbb{N}$ and $ t_1\ge 2,\,\, M \ge 8$. Let the sequence $\{t_j\}$ satisfies (\ref{lac}). Let $\{S_j\}
\subset {\rm O}_d $ be a sequence of orthogonal matrices and $ G_j = t_j S_j$. Then there exists a  real vector $ \alpha =
(\alpha_1,...,\alpha_d)$ such that
$$
 ||G_j \alpha || \ge
 \frac{1}{2^{13}d(M\log M)^{1/d}},\,\,\, \forall j \in \mathbb{N}.
$$}

Of course the constants in Theorems 1,2 may be improved.

{\bf Corollary.}\,\,\,{\it There exists an effective positive constant $\Delta$ with the following property. For any complex number $\theta =
a+bi$,
  $|\theta|>1$
  there exists
complex
   $\alpha$
such that for the the distance to the nearest Gaussian integer (in sup-norm) one has
  $$
 ||\theta^j\xi || \ge
\Delta \min\left\{\frac{ (|\theta | - 1)}{\log (2+ 1/(|\theta| - 1))}
 ; 1 \right\},\,\,\,
  \forall j\in\mathbb{N}.
$$
 }

The corollary follows from the Theorem 2 for $ d = 2$ with lacunary sequence $|\theta^j|$ and the matrix sequence $S^j$,
 $ S= \frac{1}{|\theta |}\left(
\begin{array}{cc}
a&b\cr b&-a
\end{array}
\right) $. The case $ |\theta | > 1$ easily can be reduced to the case $|\theta | > 8$.

We must note that recently Dubickas  \cite{DUBIprepr}  proved the following result. Let $t_0,t_1,t_2,... $ be a sequence of non-zero complex
numbers satisfying $|t_{j+1}| \ge a |t_j|$ with real $a > 1$ . Let $\nu  i$ be a complex number. Then there exists a complex number $\alpha$
such that  the numbers $\alpha t_j, j = 0,1,2,... $ all lie outside the union of open squares centered at $\nu +\mathbb{Z}[i]$ whose sides,
parallel to real and imaginary axis are equal to $(a-1)/20$ for $a \in (1,11-4\sqrt{5}]$ and $(a-2)/(a-1)$ for $ a > 11-4\sqrt{5}$.

{\bf 4. The principal lemma}

{\bf Lemma 1.}\,\,\, {\it Let $\mu \left(\bigcap_{j\le i} A_j^c\right) \neq 0$ Then
$$\mu \left(
A_{i+h}\bigcap (\bigcap_{j\le i} A_j^c )\right) \le
\delta_1
 \mu \left(\bigcap_{j\le i} A_j^c\right).
$$
}

Proof. The proof is identical to the proof of formula (3.9) from \cite{P}. We have
$$
\mu\left(
A_{i+h}\bigcap (\bigcap_{j\le i} A_j^c)\right)=
\sum_{\nu = 1}^{T_i}
\mu (A_{i+h} \bigcap I_\nu).
$$
As $\mu \left(\cap_{j\le i} A_j^c)\right) \neq 0$ the sum here is not empty and $T_i \ge 1$.
$A_{i+h}$ can be covered by a union of closed dyadic boxes of the form
$$
\left[\frac{b_1}{2^{l_{j+h}}}, \frac{b_1+1}{2^{l_{j+h}}}\right] \times . . . \times
\left[\frac{b_d}{2^{l_{j+h}}}, \frac{b_d+1}{2^{l_{j+h}}}\right] .
$$
 Let
$J$ be a dyadic  box from this covering of $ A_{i+h}   $  (its measure is equal to $ 2^{-dl_{i+h}}$). If $ \mu (J\cap I_\nu ) \neq 0$ then $
J\subseteq I_\nu$. Let
$$
I_\nu= \left[\frac{b_1}{2^{l_{i}}}, \frac{b_1+1}{2^{l_{i}}}\right] \times . . . \times
\left[\frac{b_d}{2^{l_{i}}}, \frac{b_d+1}{2^{l_{i}}}\right] $$ and the box $J$ appears from the
covering of the box $E(i+h, {\bf a})$. Then
$$
\frac{\bf a}{t_{i+h}} \in
 \left(\frac{b_1}{2^{l_{i}}}-\frac{\delta}{t_{i+h}},
\frac{b_1+1}{2^{l_{i}}}+\frac{\delta}{t_{i+h}}\right) \times . . . \times
\left(\frac{b_d}{2^{l_{i}}}-\frac{\delta}{t_{i+h}},
\frac{b_d+1}{2^{l_{i}}}+\frac{\delta}{t_{i+h}}\right)
.
$$
Let $W_\nu$ be the number of integer points $\bf a$ satisfying the condition above. Then
\begin{equation}
W_\nu \le \left( \left\lfloor \left(\frac{1}{2^{l_i}}+\frac{2\delta}{t_{i+h}}\right) t_{i+h} \right\rfloor +1 \right)^d \le (2^{-l_i} t_{i+h} +
2)^d. \label{W} \end{equation}
 Now we deduce
\begin{equation}
 \mu\left( A_{i+h}\bigcap \left(\bigcap_{j\le i} A_j^c\right)\right)\le 2^d\sum_{\nu =1}^{T_i}
\frac{2^d W_\nu}{2^{dl_{i+h}}} \le \sum_{\nu =1}^{T_i} \frac{(2^{-l_i} t_{i+h} + 2)^d}{2^{dl_{i+h}}}
\le \label{1}
\end{equation}
$$
\le
2^{2d} \sum_{\nu =1}^{T_i} \left(\max \left(\frac{2^{-l_i} t_{i+h}}{2^{l_{i+h}}},
\frac{1}{2^{l_{i+h}-1}}\right)\right)^d
\le
2^{2d}
\left(
\sum_{\nu =1}^{T_i}
\mu (I_\nu )\left(\frac{  t_{i+h}}{2^{l_{i+h}}}\right)^d
+
\sum_{\nu
=1}^{T_i}
\left(\frac{1}{2^{l_{i+h}-1}}\right)^d
\right)
.
 $$
But
$$
\sum_{\nu =1}^{T_i} \mu (I_\nu )\left(\frac{  t_{i+h}}{2^{l_{i+h}}}\right)^d = \mu \left(\bigcap_{j\le
i} A_j^c\right) \left(\frac{  t_{i+h}}{2^{l_{i+h}}}\right)^d. $$ Applying (\ref{DD}) we have
$l_{i+h}\ge \log_2 \frac{t_{i+h}}{2\delta} -1$ and
\begin{equation}
\sum_{\nu =1}^{T_i} \mu (I_\nu )\left(\frac{  t_{i+h}}{2^{l_{i+h}}}\right)^d
\le
(4\delta)^d \mu \left(\bigcap_{j\le i} A_j^c\right).
\label{2}
\end{equation}

From another hand
$$
\sum_{\nu =1}^{T_i} \left(\frac{1}{2^{l_{i+h}-1}}\right)^d
=
\frac{2^dT_i}{2^{dl_{i+h}}}=
2^d
\mu \left(\bigcap_{j\le i} A_j^c\right)
\left(\frac{2^{l_i}}{2^{l_{i+h}}}\right)^d
$$
as $\mu \left(\bigcap_{j\le i} A_j^c\right)= T_i 2^{-dl_i}$. Now from (\ref{DD}) we have $
\frac{2^{l_i}}{2^{l_{i+h}}}\le 2\frac{t_i}{t_{i+h}}$. Applying (\ref{HH}) we deduce
\begin{equation}
\sum_{\nu\left(\bigcap_{j\le i} A_j^c\right)
=1}^{T_i}
\left(\frac{1}{2^{l_{i+h}-1}}\right)^d
\le
2^{d+1}\delta^d \mu \left(\bigcap_{j\le i} A_j^c\right).
\label{3}
\end{equation}
Lemma 1 follows from (\ref{1},\ref{2},\ref{3}).

{\bf 5.  Proof of Theorem 1.}

The arguments of this section are te same as the arguments  from the variant of the Lovasz locall  lemma used in \cite{P}.

We shall prove by induction that the inequality
\begin{equation}
\mu
\left(\bigcap_{j\le i} A_j^c\right)\ge
\frac{1}{2}
\mu
\left(\bigcap_{j\le i-h} A_j^c\right)> 0
,\,\,\, \label{I}
\end{equation}
holds for all natural $i$.

 1. The base of induction.
For  $ i\le 0$ we define $ A_i^c = [0,1]^d$. Then the statement is trivial for $i\le 0$.
  We shall check
(\ref{I}) for $ 0\le i \le h$. It is sufficient to see that
\begin{equation}
\mu \left(
\bigcap_{1\le
j\le h}
A_j^c
\right)\ge \frac{1}{2}.
\label{bs}
\end{equation}
But
$$
\mu \left(
\bigcap_{1\le
j\le h}
A_j^c
\right)\ge
1-
\sum_{j=1}^h
\mu (
A_j
).
$$
We must take into account (\ref{mu}). Now (\ref{bs})  follows from (\ref{de}).

2. The inductive step. We suppose (\ref{I}) to be true for all $ i \le t$.

We have
$$
 \bigcap_{j\le t+1} A_j^c=
 \left( ... \left(
\left(\bigcap_{j\le t+1- h} A_j^c\right)
\setminus A_{t+1-h+1}\right)\setminus\cdots\right)\setminus A_{t+1} .
$$
Then
\begin{equation}
\mu \left(
 \bigcap_{j\le t+1} A_j^c\right)
\ge \mu \left(\bigcap_{j\le t+1- h} A_j^c\right)\ - \sum_{ v=1}^h \mu\left( A_{ t - h + 1 + v}
\bigcap
 \left(\bigcap_{j\le t+1- h} A_j^c\right)\right)
.
\label{A}
\end{equation}
 Note that for the values of $v$ under consideration we have $ t+1-h \ge t+1+v-2h$. It means
that
$$
\bigcap_{j\le t+1 - h}
A_j^c
\subseteq \bigcap_{j\le t+1+v - 2h}
A_j^c .$$
 Hence by Lemma 1
$$
\mu\left( A_{ t - h + 1 + v}
\bigcap
 \left(\bigcap_{j\le t+1- h} A_j^c\right)
\right)
 \le
\mu\left( A_{ t - h + 1 + v}
\bigcap
 \left(\bigcap_{j\le t+1+v - 2h} A_j^c\right)
\right)
\le
$$
$$
\le
\delta_1
 \mu
 \left(\bigcap_{j\le t+1+v -2h} A_j^c\right).
$$
So for $ 1\le v \le h$ we have
$$
\mu \left( A_{t-h +1+ v}\bigcap\left( \bigcap_{ j\le t+1-h}A_j^c\right) \right) \le \delta_1 \mu \left( \bigcap_{1\le t+1-2h+v} A_j^c \right)\le
\delta_1 \mu \left( \bigcap_{1\le t+2-2h} A_j^c \right) .
$$
But $ h \ge 2$ and from our inductive hypothesis for $t+2-h$ we have
$$
\mu\left(
\bigcap_{j\le
t+1-h}
A_j^c\right)
\ge
\mu\left(
\bigcap_{j\le
t+2-h}
A_j^c\right)
\ge \frac{1}{2}
\mu\left(
\bigcap_{j\le
t+2-2h}
A_j^c\right)
$$
and hence
$$
\mu\left(
\bigcap_{j\le
t+2-2h}
A_j^c\right)
\le 2
\mu\left(
\bigcap_{j\le
t+1-h}
A_j^c\right).
$$
Now
\begin{equation}
\mu \left( A_{t+1+h - v}\bigcap\left( \bigcap_{ j\le t+1-h}A_j^c\right) \right) \le 2{\delta_1} \mu \left( \bigcap_{1\le t+1-h} A_j^c \right) .
\label{B}
\end{equation}
So from (\ref{A},\ref{B}) we have
$$
\mu \left(
 \bigcap_{j\le t+1} A_j^c\right)
\ge
\left(1-2{\delta_1h}\right)
\mu \left(\bigcap_{j\le t+1- h} A_j^c\right)\
$$
and (\ref{I}) follows from the inequality (\ref{HD}). Now from (\ref{I}) and the compactness of the
sets $A_i^c$ Theorem 1 follows.

{\bf 6. Comments on the proof of Theorem 2.}

The proof of the first inequality from  theorem 2 is quite similar to the proof of the Theorem 1. Instead of parameters $\delta, \delta_1, h$ we
should use
$$
\delta' = \frac{1}{2^{13}d(M\log M)^{1/d}} ,\,\,\,  \delta_1' =
2^{4d+1}d^d(\delta ')^d,\,\,\,
 h' =\lceil 2^{3d+3} M\log M\rceil.
$$
The inequalities analogical to (\ref{de},\ref{HD},\ref{HH}) are satisfied.  We must deal with the sets $ E'(j,{\bf a}) = S_j^{-1} E(j,{\bf a})$
and their covering by dyadic boxes $A_j'$. Then instead of (\ref{mu}) we get
$$
\mu (A_j') \le \left( \frac{2(\lceil\sqrt{d}\rceil+1)\delta '}{t_j}\right)^d (\lceil \sqrt{d}t_j\rceil +1)^d\le (16\delta ')^d .
$$
 instead of the Lemma 1 we get

{\bf Lemma 2.}\,\,\, {\it Let $\mu \left(\bigcap_{j\le i} (A_j')^c\right) \neq 0$ Then
$$\mu \left(
A_{i+h}'\bigcap (\bigcap_{j\le i} (A_j')^c )\right) \le \delta_1'
 \mu \left(\bigcap_{j\le i} (A_j')^c\right).
$$
}

The sketch of the proof is as follows. For the number $W'_\nu$ of the boxes $  E'(j+h,{\bf a})$ intersecting the box $I_\nu$ we have
$$
W_\nu ' \le (2^{-l_i} t_{i+h} + 2\sqrt{d})^d
$$
instead of (\ref{W}). Then instead of (\ref{1}) we get
$$
 \mu\left( A_{i+h}'\bigcap \left(\bigcap_{j\le i} (A_j')^c\right)\right)\le (2\sqrt{d})^d\sum_{\nu =1}^{T_i}
\frac{2^d W_\nu'}{2^{dl_{i+h}}} \le
$$
$$
\le
 2^{2d} d^d \left( \sum_{\nu =1}^{T_i} \mu (I_\nu )\left(\frac{  t_{i+h}}{2^{l_{i+h}}}\right)^d + \sum_{\nu
=1}^{T_i} \left(\frac{1}{2^{l_{i+h}-1}}\right)^d \right) .
 $$
 From these inequalities we deduce Lemma 2.
 Now the proof of the Theorem 2 follows section 5 word by word.

\vskip+2.0cm

author: Nikolay Moshchevitin

\vskip+0.5cm

e-mail: moshchevitin@mech.math.msu.su, moshchevitin@rambler.ru

\end{document}